\documentclass[12pt]{article}
\usepackage{amsfonts}
\usepackage{mathrsfs}

\usepackage{amsmath,amssymb,amsfonts,theorem,makeidx,latexsym,epsfig,color,subfigure,algorithm,algorithmic}\usepackage[numbers,sort&compress]{natbib}

\newtheorem{defn}{Definition}[section]

\newtheorem{lem}[defn]{Lemma}

{\theorembodyfont{\rmfamily}

}

\newtheorem{thm}[defn]{Theorem}

\newtheorem{prop}[defn]{Proposition}

\numberwithin{equation}{section}

\def\bp{{\noindent\bf Proof. \ }}

\def\ep{\hfill$\square$\par\bigskip}

\title{Disjunctive domination in graphs with minimum degree at least two \footnote{The research is supported
by NSFC (No.11301440) and Natural Science Foundation of Fujian
Province (CN)(2015J05017)}}

\author{{\large Wei Zhuang\thanks{Corresponding author; E-mail: zhuangweixmu@163.com}} \\
{\it \normalsize School of Applied Mathematics},
\\{\it \normalsize  Xiamen University of Technology, Xiamen 361024, P.R.China }}

\date{}

\begin{document}

\maketitle

\begin{abstract} A set $D$ of vertices in $G$ is a disjunctive dominating set in $G$ if
every vertex not in $D$ is adjacent to a vertex of $D$ or has at
least two vertices in $D$ at distance $2$ from it in $G$. The
disjunctive domination number, $\gamma^{d}_2(G)$, of $G$ is the
minimum cardinality of a disjunctive dominating set in $G$. In this
paper, we show that if $G$ be a graph of order at least $3$,
$\delta(G)\geq 2$ and with no component isomorphic to any of eight
forbidden graphs, then $\gamma^{d}_2(G)\leq \frac{|G|}{3}$.
Moreover, we provide an infinite family of graphs attaining this
bound. In addition, we also study the case that $G$ is a claw-free
graph with minimum degree at least two.
\end{abstract}

\begin{minipage}{150mm}



\end{minipage}

\section{Introduction}
Let $G=(V, E)$ be a simple graph, and $v$ be a vertex in $G$. The
\emph{open neighborhood} of $v$ is $N(v)=\{u\in V|uv\in E\}$ and the
\emph{closed neighborhood} of $v$ is $N[v]=N(v)\cup \{v\}$. The
\emph{degree} of a vertex $v$ is $d(v)=|N(v)|$. A \emph{leaf} of $G$
is a vertex of degree $1$ and a \emph{support vertex} of $G$ is a
vertex adjacent to a leaf. A \emph{linkage} in $G$ is a path such
that each of its internal vertices has degree $2$ in $G$. A
\emph{subdivided star} $K_{1, t}^{*}$ is a tree obtained from a star
$K_{1, t}$ on at least two vertices by subdividing each edge exactly
once. Let $C_{s, t}$ be the graph obtained from a cycle $C_s$ and a
path $P_t$ by identifying a leaf of $P_t$ with some vertex of $C_s$,
where $s\geq 3$ and $t\geq 2$.

A \emph{dominating set} in a graph $G$ is a set $S$ of vertices of
$G$ such that every vertex in $V(G)\setminus S$ is adjacent to at
least one vertex in $S$. The domination number of $G$, denoted by
$\gamma(G)$, is the minimum cardinality of a dominating set of $G$.
The literature on the subject of domination parameters in graphs up
to the year 1997 has been surveyed and detailed in the two books
\cite{Haynes1, Haynes2}.

Motivated by the concepts of distance domination and exponential
domination (see, \cite{Anderson, Dankelmann, Henning0}), Goddard,
Henning and McPillan \cite{Goddard} introduced and studied the
concept of disjunctive domination in a graph. A set $S$ of vertices
in a graph $G$ is a \emph{disjunctive dominating set}, abbreviated
\emph{$2DD$-set}, in $G$ if every vertex not in $S$ is adjacent to a
vertex of $S$ or has at least two vertices in $S$ at distance $2$
from it in $G$. We say a vertex $v$ in $G$ is \emph{disjunctively
dominated}, abbreviated \emph{$2D$-dominated}, by the set $S$, if
$N[v]\cap S\neq \emptyset$ or there exist at least two vertices in
$S$ at distance $2$ from $v$ in $G$. The \emph{disjunctive
domination number} of $G$, denoted by $\gamma^{d}_2(G)$, is the
minimum cardinality of a $2DD$-set in $G$. A disjunctive dominating
set of $G$ of cardinality $\gamma^{d}_2(G)$ is called a
$\gamma^{d}_2(G)$-set. If the graph $G$ is clear from the context,
we simply write $\gamma^{d}_2$-set rather than
$\gamma^{d}_2(G)$-set. Every dominating set is a $2DD$-set. The
concept of disjunctive domination in graphs has been studied in
\cite{Goddard, Henning1, Henning2, Henning3, Panda, Jamil, Henning4}
and elsewhere.

In \cite{Goddard}, Goddard et al. proved the following theorem:

\begin{thm}[\cite{Goddard}]
If $G$ is a connected graph with $n\geq 5$, then
$\gamma^{d}_2(G)\leq \frac{n-1}{2}$.
\end{thm}

Moreover, they improved this bound when restrict the connected graph
$G$ to be a claw-free grpah.

\begin{thm}[\cite{Goddard}]
If $G$ is a connected claw-free graph of order $n$, then
$\gamma^{d}_2(G)\leq \frac{2n}{5}$ unless $G\in \{K_1, P_2, P_4,
C_4, H_3\}$, where $H_3$ is the graph obtained from $K_{1, 3}^{*}$
by adding an edge joining two of these support vertices.
\end{thm}

Our aim in this paper is to further improve the upper bound when we
restrict the graph $G$ to be a graph with minimum degree at least
two and a claw-free graph with minimum degree at least two,
respectively.

\section{Main results}

We first present the following lemmas, which are helpful for our
investgation.

\begin{prop}[\cite{Goddard}]
For $n\geq 3$, $\gamma^{d}_2(C_n)=2$ when $n=4$, and
$\gamma^{d}_2(C_n)=\lceil \frac{n}{4}\rceil$ when $n\neq 4$.
\end{prop}

Let $\mathscr{A}=\{C_4, C_5\}\cup \{G_i \mid i=1, 2, 3, 4, 5, 6\}$,
where each $G_i$ is showed in Fig.1.

\begin{center}
  \includegraphics[width=3.8in]{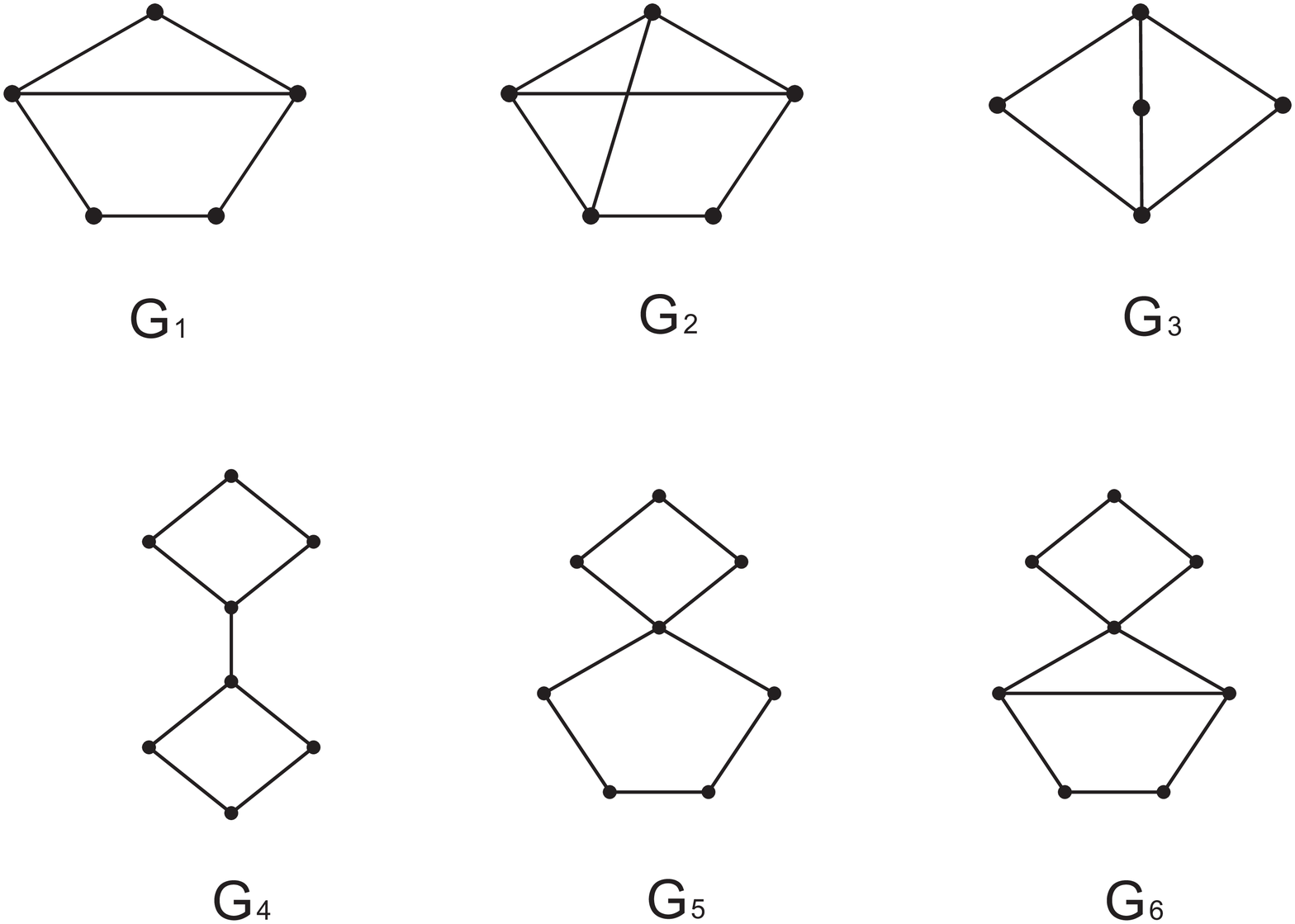}
 \end{center}
\qquad \qquad \qquad \qquad \qquad \qquad \qquad \qquad \qquad
{\small Fig.1}$\\$

\begin{lem}
Let $G$ be a connected graph of order $n\leq 8$ with $\delta(G)\geq
2$ and $G\not \in \mathscr{A}$, then $\gamma^{d}_2(G)\leq
\frac{n}{3}$.
\end{lem}

\bp This is clearly true for $n\leq 6$, so we mainly consider the
cases of $n=7$ and $8$. Let $a$ be a vertex of maximum degree in
$G$. If $d(a)=2$, then $G$ is a cycle, and the result holds. So let
$d(a)\geq 3$. If $|V(G)\setminus N[a]|=0$ or $1$, then
$\gamma^{d}_2(G)\leq 2<\frac{n}{3}$. So assume that there are two or
more vertices in $V(G)\setminus N[a]$, say $b_1, b_2, \cdots, b_s$.
Let $G'=G[\{b_1, \cdots, b_s\}]$.

Assume that $n=7$. If $s=2$ and $b_1b_2\in E(G)$, or $s=3$ and $G'$
contains at least two edges, then the result is true. Hence, we
distinguish three cases as follows.

{\flushleft\textbf{Case 1.1.}}\quad $s=2$ and $b_1b_2\not \in E(G)$.

There exist two distinct vertices $a_1'$ and $a_2'$ in $N(a)$ such
that $a_i'\in N(b_i)$, $i=1, 2$. Note that $V(G)$ can be
$2D$-dominated by $\{a_1', a_2'\}$.

{\flushleft\textbf{Case 1.2.}}\quad $s=3$ and $G'$ is independent.

In this case, $|N(a)|=3$. Since each $b_i$ has at least two
neighbors in $N(a)$, there are always two vertices in $N(a)$ which
form a $2DD$-set of $G$.

{\flushleft\textbf{Case 1.3.}}\quad $s=3$ and there is exactly one
edge in $G'$, say $b_1b_2$.

Let $a_1$ be a neighbor of $b_1$ in $N(a)$. Since $b_3$ has at least
two neighbors in $N(a)$, there exists a vertex other than $a_1$ in
$N(a)\cap N(b_3)$, say $a_2$. Take a set $D=\{a_1, a_2\}$ when
$a_2b_2\in E(G)$, and $D=\{a_1, a_3\}$ when $a_2b_2\not \in E(G)$,
where $a_3$ is the remaining vertex in $N(a)$. Clearly, $D$ is a
$2DD$-set of $G$.

The remaining case is $n=8$. It is easy to verify that if $s=2$, or
$s=3, 4$ and no edge in $G'$, or $s=3$ and there are at least two
edges in $G'$, or $s=4$ and there are at least three edges in $G'$,
we can always obtain a $2DD$-set of $G$ with cardinality $2$.

Suppose that $s=3$ and there is exactly one edge in $G'$, say
$b_1b_2$. Let $a_1, a_2$ are two neighbors of $b_3$ in $N(a)$. If
$b_1$ and $b_2$ have a common neighbor in $N(a)$, then this vertex
and $a_1$ (or $a_2$) form a $2DD$-set of $G$. Moreover, if one of
$b_1$ and $b_2$, say $b_1$, has a common neighbor with $b_3$, say
$a_1$, then $a_1$ and some vertex belonging to $N(a)\cap N(b_2)$
form a $2DD$-set of $G$. This implies that all of $b_1, b_2, b_3$
have degree two. In paricular, $b_1$ and $b_2$ are adjacent to two
distinct vertices other than $a_1$ and $a_2$, say $a_3, a_4$. It is
easy to verify that $\gamma^{d}_2(G)\leq 2$ unless $N(a)$ is a
independent set, or $a_3a_4$ is the unique edge in $G[N(a)]$. In
this two cases, $G$ is isomorphic to $G_5$ and $G_6$, respectively.

Finally, we consider the case of $s=4, 1\leq |E(G')|\leq 2$. We
distinguish three cases as follows.

{\flushleft\textbf{Case 2.1.}}\quad There is exactly one edge in
$G'$, say $b_1b_2$.

If $b_1$ and $b_2$ have a common neighbor in $N(a)$, say $a_1$, it
follows from $d(a)=3$ and $\delta(G)\geq 2$ that $a_1$ and any
vertex of $N(a)\setminus \{a_1\}$ form a $2DD$-set of $G$. So $b_1$
and $b_2$ have no common neighbor in $N(a)$, and let $a_1', a_2'$ be
the neighbors of $b_1, b_2$ in $N(a)$, respectively. Clearly, $V(G)$
is $2D$-dominated by $\{a_1', a_2'\}$.

{\flushleft\textbf{Case 2.2.}}\quad There are two nonadjacent edges
in $G'$, say $b_1b_2$ and $b_3b_4$.

Suppose that $b_1$ and $b_2$ have no common neighbor in $N(a)$, and
let $a_1\in N(b_1)\cap N(a), a_2\in N(b_2)\cap N(a)$. If
$(N(b_3)\cup N(b_4))\cap \{a_1, a_2\}\neq \emptyset$, say $b_3a_2\in
E(G)$, then either $\{b_1, b_4\}$ or $\{a_1, a_2\}$ is a $2DD$-set
of $G$. Otherwise, $(N(b_3)\cup N(b_4))\cap \{a_1, a_2\}=\emptyset$,
and then $b_3a_3, b_4a_3\in E(G)$, where $a_3$ is the remaining
vertex in $N(a)$. It means that $V(G)$ is $2D$-dominated by $\{b_1,
a_3\}$. Hence, $b_1$ and $b_2$ have a common neighbor in $N(a)$, say
$a_1'$. Similar, $b_3$ and $b_4$ have a common neighbor in $N(a)$,
say $a_2'$. Let $D=\{a_1', a_2'\}$ when $a_1'\neq a_2'$, and
$D=\{a_1', a\}$ when $a_1'=a_2'$. Clearly, $V(G)$ is $2D$-dominated
by $D$.

{\flushleft\textbf{Case 2.3.}}\quad There are two adjacent edges in
$G'$, say $b_1b_2$ and $b_2b_3$.

If there exist two distinct vertices $x$ and $y$ such that $x\in
N(a)\cap N(b_1)$ and $y\in N(a)\cap N(b_3)$, then $V(G)$ is
$2D$-dominated by$\{x, y\}$. This implies that $d(b_1)=d(b_3)=2$ and
they have a common neighbor in $N(a)$, say $a_1$. Assume that $b_4$
is adjacent to $a_1$, we have that $\{a, b_1\}$ is a $2DD$-set of
$G$. It concludes that $a_1b_4\not \in E(G)$ and $b_4$ is adjacent
to the remaining two vertices in $N(a)$. It is easy to see that
$\gamma^{d}_2(G)=2$ when $\{b_2\}\cup N(a)$ is not a independent
set, and $\gamma^{d}_2(G)=3$ when $\{b_2\}\cup N(a)$ is a
independent set. In the latter case, $G$ is isomorphic to $G_4$. \ep

\begin{lem}
Let $G$ be a connected graph, $uv\in E(G)$ and $G'$ be the graph
obtained from $G$ by subdividing $uv$ three times. If
$\gamma^{d}_2(G)\leq \frac{|G|}{3}$, then $\gamma^{d}_2(G')\leq
\frac{|G'|}{3}$.
\end{lem}

\bp Assume that $G'$ obtained from $G$ by subdividing $uv$ with
vertices $x_1, x_2, x_3$, and $S$ is a $\gamma^{d}_2$-set of $G$.
Without loss of generality, we only need to consider the following
cases.

{\flushleft\textbf{Case~1.}}\quad $u\in S$, or $u\not \in S$ and
$(N_G(u)\setminus \{v\})\cap S \neq \emptyset$. In these cases, let
$S'=S\cup \{x_3\}$ when $v\not \in S$, and $S'=S\cup \{x_1\}$ when
$v\in S$.

{\flushleft\textbf{Case~2.}}\quad $N_G[u]\cap S= \emptyset$. It
means that there are at least two vertices in $S$ at distance $2$
from $u$, then we let $S'=S\cup \{x_2\}$.

In either case, $S'$ is a $2DD$-set of $G'$. Hence,
$\gamma^{d}_2(G')\leq \gamma^{d}_2(G)+1\leq
\frac{|G|}{3}+1=\frac{|G|+3}{3}=\frac{|G'|}{3}$. \ep

Set $\mathscr{D}_1=\{C_{3, 1}, C_{3, 2}, C_{3, 3}, C_{4, 1}, C_{4,
2}, C_{4, 3}, C_{5, 1}, C_{5, 2}, C_{5, 3}\}$ and
$\mathscr{D}_2=\{C_3, C_4, C_5\}$. First, take a single vertex $v$,
$t$ graphs $A_1, A_2, \cdots, A_t\in \mathscr{D}_1$ and $s$ graphs
$B_1, B_2, \cdots, B_s\in \mathscr{D}_2$, where $s, t\geq 0$ and
$s+t\geq 2$. Let $G$ be a graph obtained by identifying the leaf of
each $A_i$ and a vertex of each $B_i$ with the vertex $v$, let
$\mathscr{F}$ be a family consisting of all such graphs $G$. In
particular, we call the each of these $A_i$ and $B_i$ the
\emph{gadget} of $G$. We give an example as follows.

\begin{center}
  \includegraphics[width=6in]{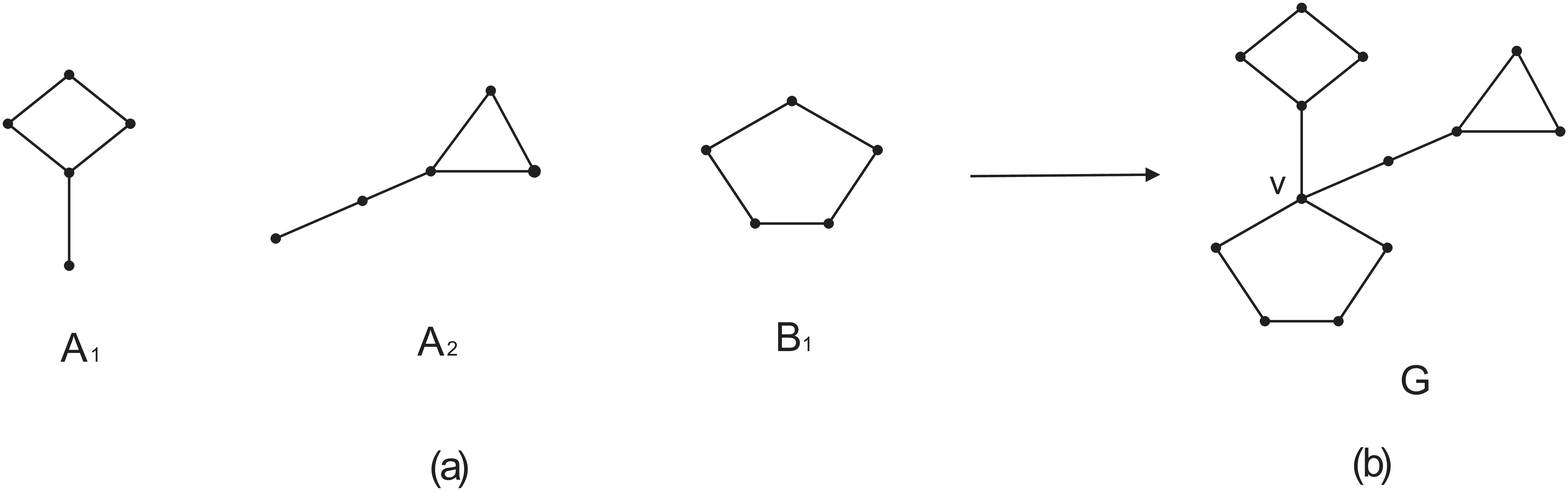}
 \end{center}
\qquad \qquad \qquad \qquad \qquad \qquad \qquad \qquad \qquad
{\small Fig.2} $\\$

\begin{lem}
If a graph $G\in \mathscr{F}\setminus \{G_4, G_5\}$ $(G_4, G_5$ are
the two graphs which are shown in Fig.1$)$, then
$\gamma^{d}_2(G)\leq \frac{|G|}{3}$.
\end{lem}

\bp The proof is by induction on $|V(G)|$. If $G$ contains exactly
two gadgets, then it is easy to verify that the result holds. So we
consider the case that $G$ contains at least three gadgets.

Assume that $G$ contains a $4$-cycle $C$ as its gadget, we removing
all $2$-degree vertices of this cycle, denote the resulting graph by
$G'$. If $G'$ is not isomorphic to $G_4$ or $G_5$, then by
induction, $\gamma^{d}_2(G')\leq \frac{|G'|}{3}$. Let $S$ be a
$\gamma^{d}_2$-set of $G'$. Then, $S\cup \{u\}$ is a $2DD$-set of
$G$, where $u$ is the vertex at distance $2$ from $v$ in $C$. Hence,
we have that $\gamma^{d}_2(G)\leq \gamma^{d}_2(G')+1\leq
\frac{|G'|}{3}+1=\frac{|G|}{3}$. On the other hand, if $G'$ is
isomorphic to $G_4$ or $G_5$, then it is easy to see that
$\gamma^{d}_2(G)\leq \frac{|G|}{3}$. Similar, if there exists a
gadget of $G$ which is isomorphic to one of $\{C_{3, 1}, C_{4, 3},
C_{5, 2}\}$, we are done. Next, suppose that $T_1$ and $T_2$ are two
of the gadgets of $G$.

If one of $T_1$ and $T_2$, say $T_1$, belongs to $\{C_{4, 2}, C_{5,
1}\}$, then removing all vertices of $T_2$ except $v$. Denote the
resulting graph by $G'$. Due to the existence of $T_1$, there must
be a $\gamma^{d}_2$-set of $G'$, say $S_1$, such that $|S_1|\leq
\frac{|G'|}{3}$ and $v\in S_1$. Based on this, the set $S$ can
always be extended to a $2DD$-set of $G$ with cardinality at most
$\frac{|G|}{3}$.

Hence, we assume that each gadget of $G$ belongs to $\{C_3, C_5,
C_{3, 2}, C_{3, 3}, C_{4, 1}, C_{5, 3}\}$. Moving both of $T_1$ and
$T_2$ except $v$, after that, adding a $3$-cycle $C_1$ and
identifying a vertex of this cycle with $v$. Denote the resulting
graph by $G''$. Let $S_2$ be a $\gamma^{d}_2$-set of $G''$. Note
that $N[v]\cap S_2\neq \emptyset$. Let $S'=(S_2\setminus V(C_1))\cup
\{v\}$ when $S_2\cap V(C_1)\neq \emptyset$, and $S'=S_2$ when
$S_2\cap V(C_1)=\emptyset$, it is easy to see that $S'$ can always
be extended to a $2DD$-set of $G$ with cardinality at most
$\frac{|G|}{3}$. \ep

Let $H$ be a graph satisfying its all vertices have degree at least
two except possibly one vertex $v$. We give the definitions of some
graphs as follows, the related graphs are shown in Fig.3.

\begin{center}
  \includegraphics[width=6.2in]{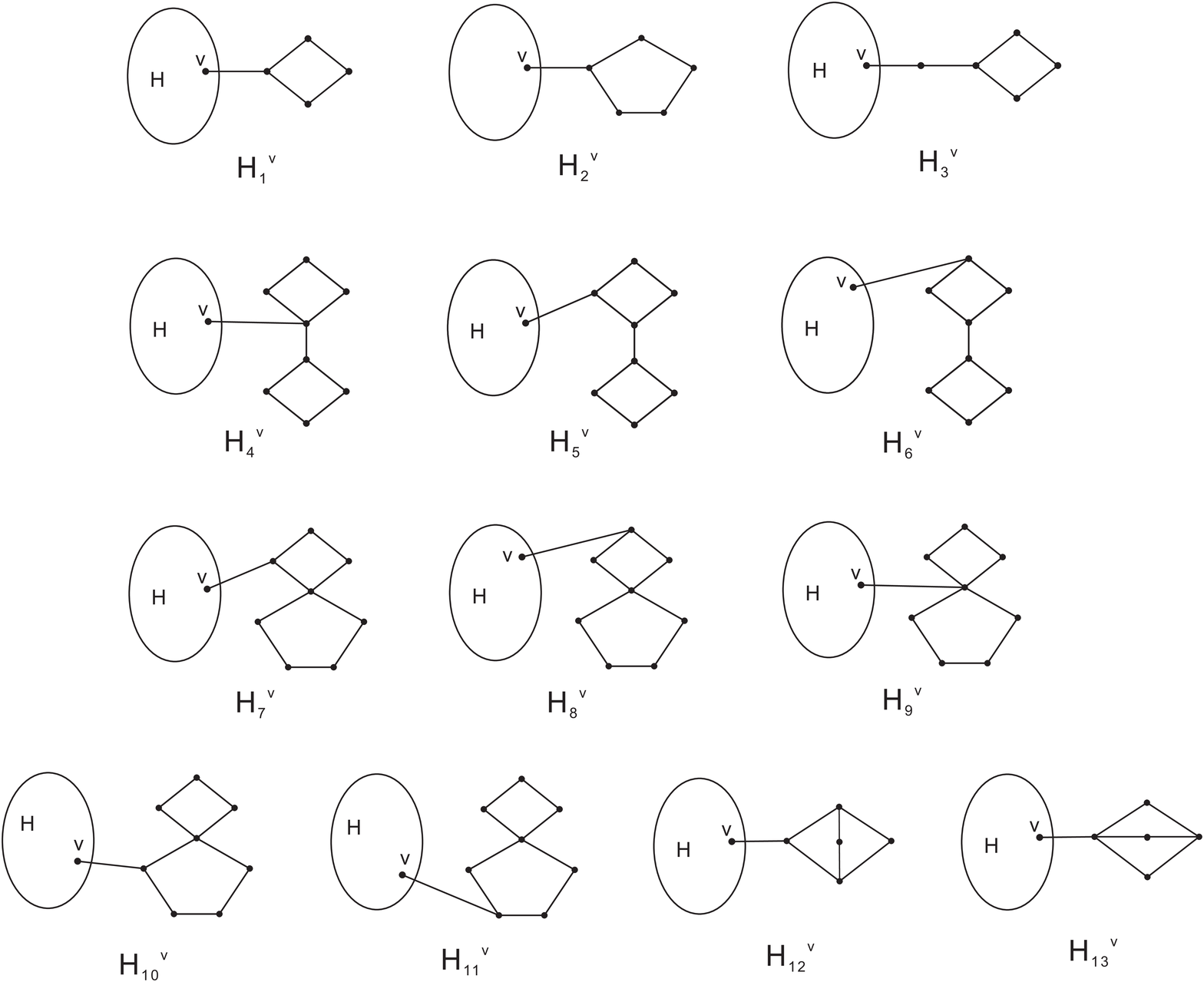}
 \end{center}
\qquad \qquad \qquad \qquad \qquad \qquad \qquad \qquad \qquad
{\small Fig.3} $\\$

$\bullet$ $H_1^{v}$ (respectively, $H_2^{v}$) denotes the graph
obtained from $H$ by adding a $C_4$ (respectively, $C_5$), and
joining $v$ to a vertex of $C_4$ (respectively, $C_5$).

$\bullet$ $H_3^{v}$ denotes the graph obtained from $H$ by adding a
$C_{4, 1}$, and joining $v$ to the leaf of $C_{4, 1}$.

$\bullet$ $H_i^{v}$ ($4\leq i\leq 6$) denotes the graph obtained
from $H$ by adding a $G_4$, and joining $v$ to some vertex of $G_4$.

$\bullet$ $H_j^{v}$ ($7\leq j\leq 11$) denotes the graph obtained
from $H$ by adding a $G_5$, and joining $v$ to some vertex of $G_5$.

$\bullet$ $H_j^{v}$ ($12\leq j\leq 13$) denotes the graph obtained
from $H$ by adding a $G_3$, and joining $v$ to some vertex of $G_3$.

\begin{lem}
$(i)$ If $\gamma^{d}_2(H_1^{v})\leq \frac{|H_1^{v}|}{3}$, then
$\gamma^{d}_2(H_i^{v})\leq \frac{|H_i^{v}|}{3}$, where $5\leq i\leq
8$.

$(ii)$ If $\gamma^{d}_2(H_2^{v})\leq \frac{|H_2^{v}|}{3}$, then
$\gamma^{d}_2(H_j^{v})\leq \frac{|H_j^{v}|}{3}$, where $9\leq j\leq
13$.

$(iii)$ If $\gamma^{d}_2(H_3^{v})\leq \frac{|H_3^{v}|}{3}$, then
$\gamma^{d}_2(H_4^{v})\leq \frac{|H_4^{v}|}{3}$.
\end{lem}

\bp We only need to consider the first case of $(i)$, and the proofs
of the other cases (including $(ii)$ and $(iii)$) are similar.

Denote the vertices of $V(H_5^{v}-H)$ by $\{x_1, x_2, \cdots, x_8\}$
(see Fig.4). Note that $H_1^{v}$ is obtained from $H_5^{v}$ by
deleting $x_5, x_6, x_7$ and $x_8$. Let $S$ be a $\gamma^{d}_2$-set
of $H_1^{v}$. Without loss of generality, we only need to consider
the following cases.

{\flushleft\textbf{Case~1.}}\quad $x_4\in S$, or $x_4\not \in S$ and
$N(x_4)\cap S \neq \emptyset$ in $H_1^{v}$. In these cases, let
$S'=S\cup \{x_7\}$.

{\flushleft\textbf{Case~2.}}\quad $N[x_4]\cap S=\emptyset$ in
$H_1^{v}$. It means that there are at least two vertices in $S$ at
distance $2$ from $x_4$ in $H_1^{v}$, Since $vx_1$ is a cut edge of
$H_1^{v}$, $\{v, x_2\}\subseteq S$. Let $S'=(S\cup \{x_2\})\cup
\{x_3, x_7\}$.

In either case, $S'$ is a $2DD$-set of $H_5^{v}$. Hence,
$\gamma^{d}_2(H_5^{v})\leq \gamma^{d}_2(H_1^{v})+1\leq
\frac{|H_1^{v}|}{3}+1=\frac{|H_1^{v}|+3}{3}<\frac{|H_5^{v}|}{3}$.
\ep

\begin{center}
  \includegraphics[width=4in]{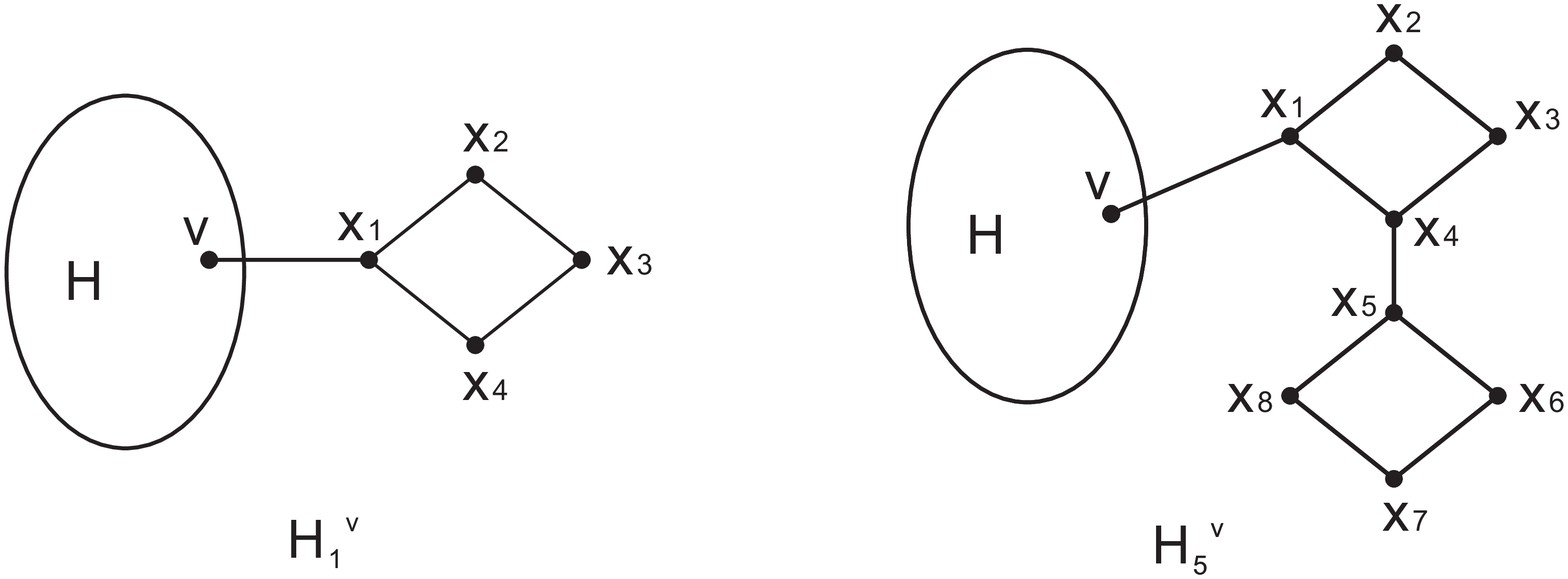}
 \end{center}
\qquad \qquad \qquad \qquad \qquad \qquad \qquad \qquad \qquad
{\small Fig.4} $\\$

Let $H$ be a graph with minimum degree at least two and $v\in V(H)$,
the graph $G$ obtained from $H$ and a graph $H'\in \mathscr{A}$ by
joining $v$ to some vertex of $H'$, say $u$, and by subdividing the
edge $uv$ $k$ times, where $k\geq 0$. The subgraph of $G$ induced by
$(V(G)\setminus V(H))\cup \{v\}$ is called a \emph{special pendant
subgraph} of $G$ attached at the vertex $v$, or simply, a
\emph{special pendant subgraph} of $v$. In particular, each vertex
of $V(G)\setminus V(H)$ is called an \emph{internal vertex} of the
special pendant subgraph. If $C_m$ is an induced cycle of a graph
$G$ such that one of its vertices, say $v$, is a cut vertex of $G$,
and $d_G(v)\geq 4$, then this $C_m$ is called a \emph{pendant cycle}
of $G$ attached at the vertex $v$, or simply, a \emph{pendant cycle}
of $v$. We take an example to make it easier for reader to
understand these definitions. In Fig.5, the vertex $u$ has a special
pendant subgraph which is isomorphic to a $C_{4, 2}$, and the vertex
$v$ has a pendant cycle which is isomorphic to a $C_3$.

\begin{center}
  \includegraphics[width=2.8in]{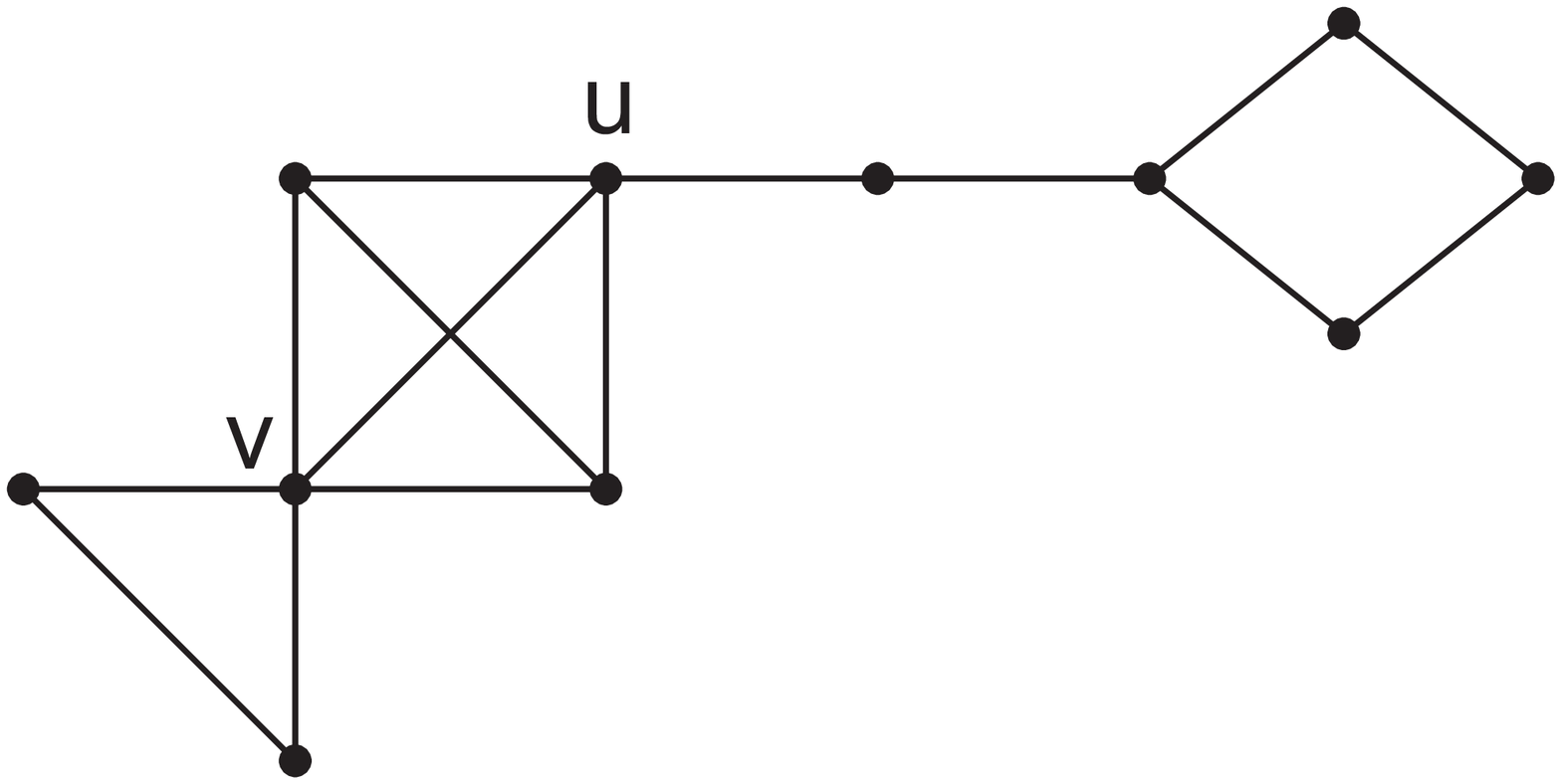}
 \end{center}
\qquad \qquad \qquad \qquad \qquad \qquad \qquad \qquad \qquad
{\small Fig.5} $\\$

\begin{thm}
Let $G$ be a graph of order at least $3$, $\delta(G)\geq 2$ and with
no component isomorphic to a graph belonging to $\mathscr{A}$,
then $\gamma^{d}_2(G)\leq \frac{|G|}{3}$.
\end{thm}

\bp We use induction on $|V(G)|+|E(G)|$. If $|V(G)|\leq 8$, then by
Lemma~2.2, the result holds. So let $|E(G)|\geq |V(G)|\geq 9$,
assume that the result holds for every graph $G'$ with
$\delta(G')\geq 2$ satisfying the conditions that
$|V(G')|+|E(G')|<|V(G)|+|E(G)|$ and no component isomorphic to a
graph belonging to $\mathscr{A}$. Note that $\gamma^{d}_2(G)\leq
\gamma^{d}_2(G-e)$ for every edge $e\in E(G)$, so we let $|E(G)|$ be
as small as possible.

If $G$ is not a connected graph, let $G_1, G_2, \cdots, G_k$ be the
components of $G$. By induction, $\gamma^{d}_2(G_i)\leq
\frac{|G_i|}{3}$ for $1\leq i\leq k$. And then,
$\gamma^{d}_2(G)=\sum \limits_{i=1}^{k} \gamma^{d}_2(G_i)\leq \sum
\limits_{i=1}^{k} \frac{|G_i|}{3}=\frac{|G|}{3}$. Hence, we consider
the case that $G$ is connected. If $\Delta(G)=2$, then $G$ is a
cycle. Since $G\neq C_4$ or $C_5$, by Proposition~2.1, we are done.
So we consider the case of $\Delta(G)\geq 3$, and let $R=\{v\in V(G)
\mid d(v)\geq 3\}$. From the choice of $G$ and the condition of
Theorem~2.6, we have that:

$(*)$ No cycle in $G$ has a chord;

$(**)$ If there are two adjacent vertices belong to $R$, say $x_1,
x_2$, then $G-x_1x_2$ is disconnected, and that at least one of the
components is isomorphic to a graph of $\mathscr{A}$.

Moreover, combining the above conclusions with Lemma~2.3, 2.5, we
have that each special pendant subgraph of $G$ is isomorphic to some
$C_{r, k}$, where $3\leq r \leq 5, 1\leq k\leq 3$, and each pendant
cycle of $G$ is isomorphic to some $C_m$, where $3\leq m\leq 5$. In
particular, each linkage in $G$ has at most two internal vertices.

Let $A=\{v \mid v$ is an internal vertex of a special pendant
subgraph of $G$ and $v\in R\}$, $B=R-A$. If $|B|=1$, it follows from
Lemma~2.4 that the result holds. Thus, $|B|\geq 2$. Moreover, by the
choice of $G$ and the argument as above, $B$ is independent.

{\flushleft\textbf{Claim~1.}}\quad For a vertex $x\in B$, there is
at most one special pendant subgraph (or pendant cycle) attached at
it. In particular, this pendant subgraph (or pendant cycle) is not
isomorphic to any of $\{C_4, C_{3, 1}, C_{4, 3}, C_{5, 2}\}$.$\\$

If there is a pendant cycle $C$ of order $4$ attached at $x$, we
removing all $2$-degree vertices of this cycle, denote the resulting
graph by $G'$. Clearly, $x$ has degree at least two in $G'$. If $G'$
is not isomorphic to any graph of $\mathscr{A}$, then by induction,
$\gamma^{d}_2(G')\leq \frac{|G'|}{3}$. Let $S$ be a
$\gamma^{d}_2$-set of $G'$. Then $S\cup \{x'\}$ is a $2DD$-set of
$G$, where $x'$ is the vertex at distance $2$ from $x$ in $C$. If
$G'$ is isomorphic to some graph of $\mathscr{A}$, then it is easy
to verify that $\gamma^{d}_2(G)\leq \frac{|G|}{3}$. Similar, if
$C_{3, 1}$, or $C_{4, 3}$, or $C_{5, 2}$ is a special pendant
subgraph of $x$, we are done.

Next, suppose that $T_1$ and $T_2$ is the special pendant subgraphs
(or pendant cycles) of $x$. If one of $T_1$ and $T_2$, say $T_1$,
belongs to $\{C_{4, 2}, C_{5, 1}\}$, then removing all vertices of
$T_2$ except $x$. Denote the resulting graph by $G'$. It is easy to
see that $x$ has degree at least two in $G'$. Due to the existence
of $T_1$, there must be a $\gamma^{d}_2$-set of $G'$, say $S$, such
that $|S|\leq \frac{|G'|}{3}$ and $x\in S$. Based on this, the set
$S$ can always be extended to a $2DD$-set of $G$ with cardinality at
most $\frac{|G|}{3}$.

Hence, assume that $T_1, T_2\in \{C_3, C_5, C_{3, 2}, C_{3, 3},
C_{4, 1}, C_{5, 3}\}$. Moving both of $T_1$ and $T_2$ except $x$,
after that, adding a $3$-cycle $C$ and identifying a vertex of this
cycle with $x$. Denote the resulting graph by $G'$. Let $S$ be a
$\gamma^{d}_2$-set of $G'$. Note that $N_{G'}[x]\cap S\neq
\emptyset$. Let $S'=(S\setminus V(C))\cup \{x\}$ when $S\cap
V(C)\neq \emptyset$, and $S'=S$ when $S\cap V(C)=\emptyset$, it is
easy to see that $S'$ can always be extended to a $2DD$-set of $G$
with cardinality at most $\frac{|G|}{3}$. \ep

Assume that $|B|=2$. Let $B=\{x, y\}$. From Claim~1, each of $x$ and
$y$ has at most one special pendant subgraph (or pendant cycle).
Combining this fact with $d(x), d(y)\geq 3$, there are at least two
linkages between $x$ and $y$. Let $K$ be the subgraph of $G$
consisting of $x, y$ and all linkages between them. If $|K|\geq 6$,
then $\{x, y\}$ is a $2DD$-set of $K$, and it can always be extended
to a $2DD$-set of $G$ with cardinality at most $\frac{|G|}{3}$. If
$|K|\leq 5$, then one of the following conditions holds:

$(*)$ There are two or three linkages between $x$ and $y$, and each
of them has exactly one internal vertex;

$(**)$ There are two linkages between $x$ and $y$, which contain one
and two internal vertices, respectively.

In either case, by a simple verify, the result holds.

Finally, we consider the case of $|B|\geq 3$.

{\flushleft\textbf{Claim~2.}}\quad There is no path $xaybz$ in $G$
such that $x, y, z\in B$ and $a, b$ have degree two.$\\$

Suppose that there is such a path. If $x=z$, then there are at least
two linkages between $x$ and $y$ that contain exactly one internal
vertex. Let $G'$ be a graph is obtained from $G$ by contracting $x,
y$ and all linkages between $x$ and $y$ that contain exactly one
internal vertex, into a single vertex, say $z$. Let $S'$ be a
$\gamma^{d}_2$-set of $G'$. If $G'$ is not isomorphic to any graph
of $\mathscr{A}$, then by induction, $|S'|\leq \frac{|G'|}{3}$. And
$(S'\setminus \{z\})\cup \{x, y\}$, or $S'\cup \{x\}$, or $S'\cup
\{y\}$, or $S'\cup \{w\}$ is a $2DD$-set of $G$ with cardinality at
most $\frac{|G|}{3}$, where $w\in V(G)\setminus (V(G')\cup \{x,
y\})$. If $G'$ is isomorphic to a graph of $\mathscr{A}$, then it is
easy to verify that $\gamma^{d}_2(G)\leq \frac{|G|}{3}$.

Hence we assume that $x\neq z$. Let $G'$ be the graph obtained from
$G$ by removing the edges $ax, bz$ and adding the edge $ab$. If no
component of $G'$ is isomorphic to a graph of $\mathscr{A}$, by the
induction hypothesis, $G'$ has a $\gamma^{d}$-set $S$ with
cardinality at most $\frac{|G'|}{3}$. Let $S'=(S\setminus \{a, b,
y\})\cup \{y\}$ when $S\cap \{a, b, y\}\neq \emptyset$, and $S'=S$
when $S\cap \{a, b, y\}=\emptyset$. Clearly, $S'$ is a $2DD$-set of
$G$. So we suppose that one of the components of $G'$, say $G_1'$,
is isomorphic to a graph of $\mathscr{A}$. Note that the component
containing $y$ is not isomorphic to a graph of $\mathscr{A}$. On the
other hand, if $G_1'$ contains exactly one of $x$ and $z$, then this
vertex belongs to $A$, a contradiction. It implies that $G'$ have
two components, one contains the vertex $y$, and the other contains
the two vertices $x, z$. In particular, the component containing $x,
z$ is $G_1'$. Next, we construct a graph $G''$ obtained from $G$ by
removing the edges $ay$ and $by$. $G''$ have two components $G_1$
and $G_2$, where $G_1$ contains $y$, and $G_2$ contains $x, z$. Note
that $G_2\not \in \mathscr{A}$, and it is easy to verify that
$\gamma^{d}_2(G_2)\leq \frac{|G_2|}{3}$. We distinguish the cases as
follows.

{\flushleft\textbf{Case~1.1.}}\quad $G_1\in \mathscr{A}$, or $G_1\in
\{C_{3, 1}, C_{3, 2}, C_{3, 3}, C_{4, 1}, C_{4, 2}, C_{4, 3}, C_{5,
1}, C_{5, 2}, C_{5, 3}\}$.

By a simple verification, the result is true.

{\flushleft\textbf{Case~1.2.}}\quad $d_{G''}(y)\geq 2$ and $G_1\not
\in \mathscr{A}$.

By induction, $\gamma^{d}_2(G_1)\leq \frac{|G_1|}{3}$. Therefore,
$\gamma^{d}_2(G)\leq \gamma^{d}_2(G_1)+\gamma^{d}_2(G_2)\leq
\frac{|G_1|}{3}+\frac{|G_2|}{3}= \frac{|G|}{3}$.

{\flushleft\textbf{Case~1.3.}}\quad $d_{G''}(y)=1$ and $G_1\not \in
\{C_{3, 1}, C_{3, 2}, C_{3, 3}, C_{4, 1}, C_{4, 2}, C_{4, 3}, C_{5,
1}, C_{5, 2}, C_{5, 3}\}$.

There is a path $wy_1\cdots y_ty$ such that $w\in B\setminus \{x, y,
z\}$ and each $y_i$ has degree two, where $t=1$ or $2$. We removing
all vertices of $V(G_2)\cup \{a, b, y, y_1, \cdots , y_t\}$, and
denote the resulting graph by $G^{*}$. $w$ has degree at least two
in $G^{*}$. If $G^{*}\not \in \mathscr{A}$, by induction, $G^{*}$
has a $2DD$-set $S'$ of cardinality at least $\frac{|G^{*}|}{3}$.
And $S'$ can be extended to a $2DD$-set of $G$ with cardinality at
least $\frac{|G|}{3}$. On the other hand, if $G^{*}\in \mathscr{A}$,
it is easy to verify that the result holds.\ep

{\flushleft\textbf{Claim~3.}}\quad There is no a cycle
$C=a_1P_1a_2P_2\cdots a_kP_ka_1$ such that each $a_i\in B$ and each
$P_i$ is a linkage containing two internal vertices between $a_i$
and $a_{i+1}$.$\\$

Suppose that there is such a cycle $C$, let $H$ be the subgraph of
$G$ obtained from $C$ by adding all linkages between $a_i$ and
$a_j$(if any), where $1\leq i, j\leq k$, and by adding the special
pendant subgraph (or pendant cycle) attached at $a_i$ (if any),
$i=1, 2, \cdots, k$. Clearly, the set $\{a_1, a_2, \cdots, a_k\}$
can be extended to a $2DD$-set $S_1$ of $H$ with cardinality at most
$\frac{|H|}{3}$. If $H=G$, we are done. So we assume that $H\neq G$.
Let $(\bigcup \limits_{i=1}^{k} N(a_i))\cap V(G-H)=\{u_1, u_2,
\cdots, u_t\}$. Note that each $u_i$ has degree one in $G-H$. Let
$H'$ be the graph obtained from $G-H$ by adding the edges
$u_ju_{j+1}$, $j=1, 3, \cdots, t-1$ when $t$ is even, and $H'$ be
the graph obtained from $G-H$ by adding the edges $u_ju_{j+1}$,
$j=1, 3, \cdots, t-2$, and by adding a path $v_1v_2v_3$ and the
edges $v_1u_t, v_3u_t$ when $t$ is odd. If no component of $H'$ is
isomorphic to a graph of $\mathscr{A}$, then by induction, $H'$ has
a $2DD$-set $S$ of cardinality at most $\frac{|H'|}{3}$. In the
former case, $S\cup S_1$ is a $2DD$-set of $G$ with cardinality at
most $\frac{|G|}{3}$. In the latter case, it is easy to see that the
set $D=\{v_1, v_2, v_3\}\cap S$ is not empty. Let $S'=S\setminus D$
when $|D|=1$, and $S'=(S\setminus D)\cup \{u_t\}$ when $|D|\geq 2$.
Then, $S'\cup S_1$ is desired set.

Next we let $H_1, H_2, \cdots, H_p$ are the components of $H'$, and
without loss of generality, assume that some of these components,
say $H_1, H_2, \cdots, H_q$, are isomorphic to the graphs of
$\mathscr{A}$, where $q\leq p$. Note that $u_t\not \in \bigcup
\limits_{i=1}^{q} V(H_i)$ when $t$ is odd. Let $H''=G[(\bigcup
\limits_{i=1}^{q} V(H_i))\cup V(H)]$. It is easy to see that $S_1$
can be extended to a $2DD$-set $S_2$ of $H''$ with cardinality at
most $\frac{|H''|}{3}$. Similar to the argument as above, $S_2$ can
be further extended to a $2DD$-set of $G$ with cardinality at most
$\frac{|G|}{3}$. \ep

By Claim~3, we let $P=a_1P_1a_2P_2\cdots a_kP_ka_{k+1}$ be a maximal
path in $G$, such that each $a_i\in B$ and each $P_i$ is a linkage
containing two internal vertices between $a_i$ and $a_{i+1}$. From
Claim~1, we know that each of $a_1$ and $a_{k+1}$ has at most one
special pendant subgraph (or pendant cycle). Combining the choice of
$P$ with Claim~2, there exists only one vertex $a\in B$ ($a'\in B$,
respectively) such that there is a linkage between $a$ and $a_1$
($a'$ and $a_{k+1}$, respectively). Moreover, each of the two
linkages contains exactly one internal vertex, say $x$ and $x'$.
Clearly, $a\neq a'$. In particular, each of $a_1, a_{k+1}$ has
exactly one special pendant subgraph (or pendant cycle). We
distinguish the cases as follows.

{\flushleft\textbf{Case~2.1.}}\quad $a=a_{k+1}$.

It means that $a'=a_1$. Let $H$ be the subgraph of $G$ obtained from
$P$ by adding all linkages between $a_i$ and $a_j$(if any), where
$1\leq i, j\leq k+1$, and by adding the special pendant subgraph (or
pendant cycle) attached at $a_i$ (if any), $i=1, 2, \cdots, k+1$.
Note that either $\{a_1, a_2, \cdots, a_k, a_{k+1}\}$ or $\{c, a_2,
\cdots, a_k, a_{k+1}\}$ can be extended to a $2DD$-set of $H$ with
cardinality at most $\frac{|H|}{3}$, where $c$ is a neighbor of
$a_1$ belonging to the special pendant subgraph (or pendant cycle)
of $a_1$. Similar to the argument of Claim~3, we can obtain the
desired set.

{\flushleft\textbf{Case~2.2.}}\quad $a=a_i$, $i=2, 3, \cdots, k$.

As mentioned before, each of $a_1$ and $a_{k+1}$ has exactly one
special pendant subgraph (or pendant cycle). If $a'$ is some $a_i$,
$i\in \{2, 3, \cdots, k\}$, the proof is similar to that of Claim~3.
So we assume that $a'\not \in \{a_1, a_2, \cdots, a_{k+1}\}$. It
follows from Claim~2 and Claim~3 that for any $d\in B\setminus
\{a_1, a_{k+1}, a'\}$, there is at most one linkage joining $a'$ and
$d$, and this linkage (if any) has two internal vertices. Moreover,
there is no path $a_ix_1x_2a'y_1y_2a_j$ in $G$ such that $i, j\in
\{2, 3, \cdots, k\}$ and $x_1, x_2, y_1, y_2$ have degree two.

Let $H$ be the subgraph of $G$ obtained from $Px'$ by adding all
linkages between $a_i$ and $a_j$(if any), where $1\leq i, j\leq
k+1$, and by adding the special pendant subgraph (or pendant cycle)
attached at $a_i$ (if any), $i=1, 2, \cdots, k+1$. Note that $\{a_1,
a_2, \cdots, a_k, a_{k+1}\}$ can be extended to a $2DD$-set of $H$
with cardinality at most $\frac{|H|}{3}$. Let $(\bigcup
\limits_{i=1}^{k+1} N(a_i))\cap V(G-H)=\{u_1, u_2, \cdots, u_t\}$.
Note that each $u_i$ has degree one in $G-H$. Let $H'$ be the graph
obtained from $G-H$ by adding the edges $u_ju_{j+1}$ ($j=1, 3,
\cdots, t-1$) when $t$ is even, and $H'$ be the graph obtained from
$G-H$ by adding a path $v_1v_2v_3$ and the edges $v_1u_t, v_3u_t,
u_ju_{j+1}$ ($j=1, 3, \cdots, t-2$) when $t$ is odd. If no component
of $H'$ is isomorphic to a graph of $\mathscr{A}$, we are done. So
assume that some of the components of $H'$, say $H_1, H_2, \cdots,
H_q$, are isomorphic to the graphs of $\mathscr{A}$. Note that $a'$
can not be contained in any $H_i$, $i\in \{1, 2, \cdots, q\}$. In
this case, similar to the proof of Claim~3, we can always obtain the
desired set.

{\flushleft\textbf{Case~2.3.}}\quad $a, a'\not \in B\setminus \{a_1,
a_2, \cdots, a_{k+1}\}$.

Let $H$ be the subgraph of $G$ obtained from $xPx'$ by adding all
linkages between $a_i$ and $a_j$(if any), where $1\leq i, j\leq
k+1$, and by adding the special pendant subgraph (or pendant cycle)
attached at $a_i$ (if any), $i=1, 2, \cdots, k+1$. Note that $\{a_1,
a_2, \cdots, a_k, a_{k+1}\}$ can be extended to a $2DD$-set of $H$
with cardinality at most $\frac{|H|}{3}$. Let $(\bigcup
\limits_{i=1}^{k+1} N(a_i))\cap V(G-H)=\{u_1, u_2, \cdots, u_t\}$.
Note that each $u_i$ has degree one in $G-H$. Let $H'$ be the graph
obtained from $G-H$ by adding the edges $u_ju_{j+1}$ ($j=1, 3,
\cdots, t-1$) when $t$ is even, and $H'$ be the graph obtained from
$G-H$ by adding a path $v_1v_2v_3$ and the edges $v_1u_t, v_3u_t,
u_ju_{j+1}$ ($j=1, 3, \cdots, t-2$) when $t$ is odd. In either case,
similar to the proof of Case~2.2, we can always obtain the desired
set.\ep

\begin{center}
  \includegraphics[width=3.2in]{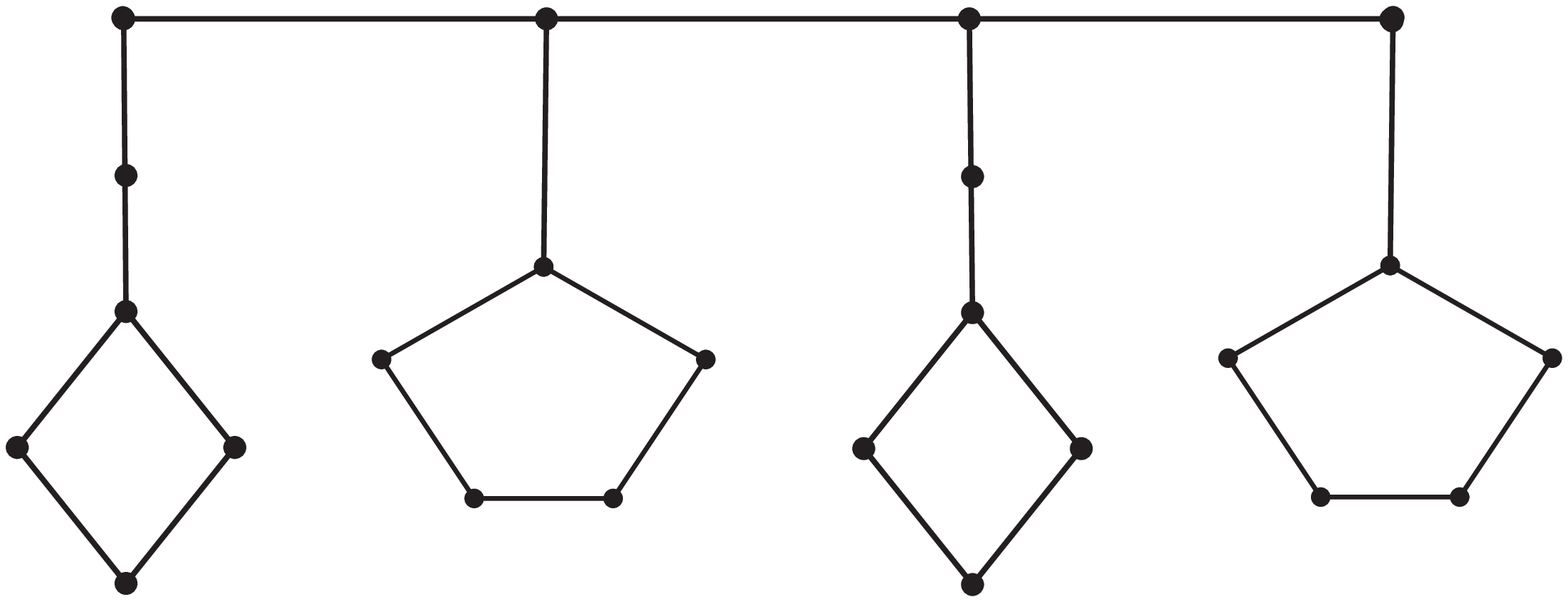}
 \end{center}
\qquad \qquad \qquad \qquad \qquad \qquad \qquad \qquad \qquad
{\small Fig.6} $\\$

Clearly, the bound of Theorem~2.6 is tight. Next, we characterize an
infinite family of graphs achieving equality for the result. Let $F$
be any connected graph of order at least $2$. For each vertex $x_i$
$(1\leq i\leq |F|)$ of $F$, identifying the leaf of a graph $T_i\in
\{C_{4, 2}, C_{5, 1}\}$ with $x_i$. Denote the resulting graph by
$G_F$, and let $\mathscr{T}$ be a family consisting of all such
graphs $G_F$. Note that $G_F$ contain $|F|$ pairwise disjoint
induced subgraphs, each of which is isomorphic to $C_{4, 2}$ or
$C_{5, 1}$. Take a $\gamma^{d}_2$-set $S$ of $G_F$, each of these
induced subgraphs contains at least two vertices belonging to $S$,
it means that $\gamma^{d}_2(G_F)\geq \frac{|G_F|}{3}$. Combining
this conclusion with Theorem~2.6, we have that
$\gamma^{d}_2(G_F)=\frac{|G_F|}{3}$.

When $F=P_4$, the graph $G_F$ is illustrated in Fig.6.

From Theorem~1.1 and Theorem~1.2, we know that the upper bound of
the disjunctive domination number of a connected graph $G$ can be
improved if $G$ is restricted to be a claw-free graph. Hence, it is
natural to consider the following question: can the upper bound of
Theorem~2.6 be improved if $G$ is restricted to be a claw-free graph
with minimum degree at least two?

\begin{center}
  \includegraphics[width=4.2in]{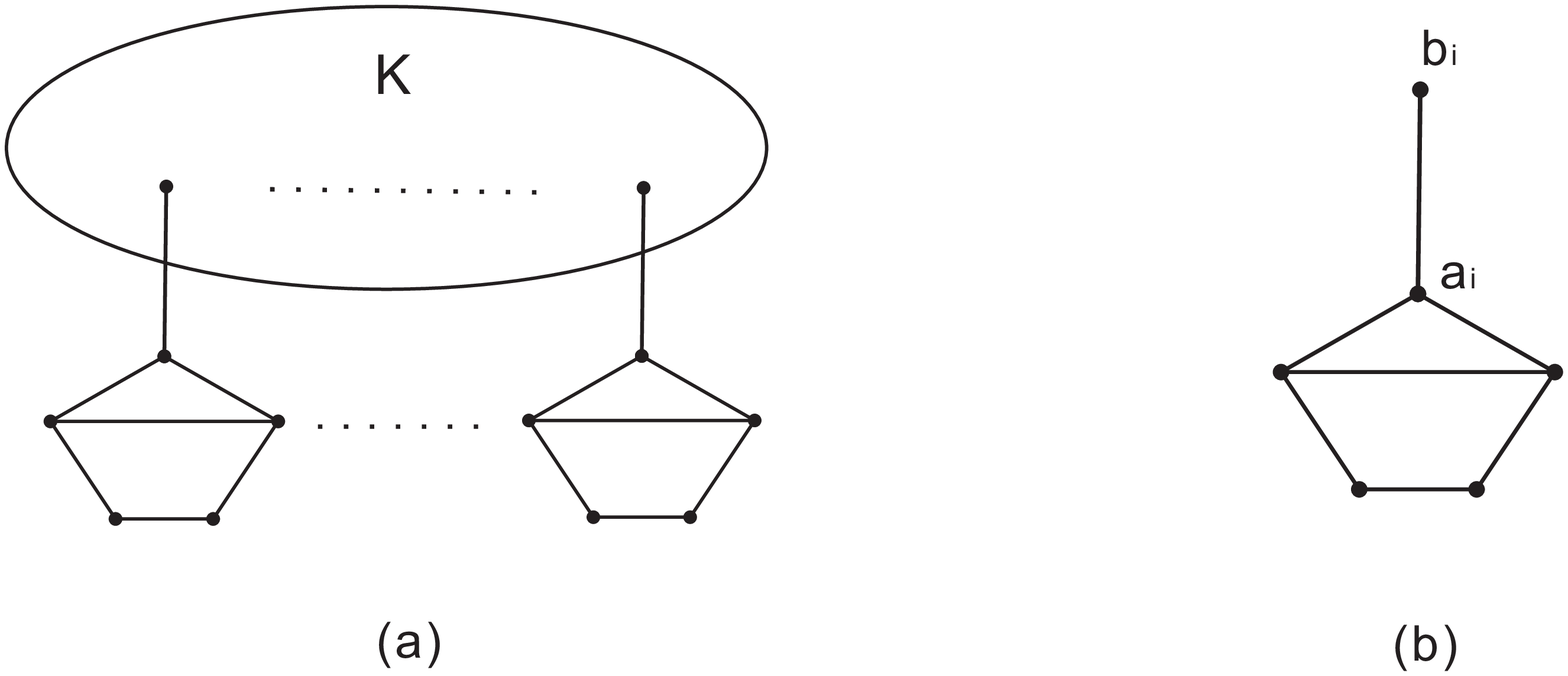}
 \end{center}
\qquad \qquad \qquad \qquad \qquad \qquad \qquad \qquad \qquad
{\small Fig.7} $\\$

Now we construct a class of graphs as follows. Let $K$ be an
complete graph of order $n\geq 2$, and let $G_K$ be the graph
obtained from $K$ by adding $n$ copies $G_{1, 1}, G_{1, 2}, \cdots,
G_{1, n}$ of $G_1$ ($G_1$ is the first graph is shown in Fig.1), by
joining $a_i$ to the $i$-th vertex $b_i$ of $K$, where $a_i$ is the
common neighbor of the two vertices of degree three in $G_{1, i}$
(see Fig.7(a)). Let $\mathscr{U}=\{G_K \mid K$ is an complete graph
of order at least two$\}$. Take a $\gamma^{d}_2$-set $S$ of $G_K$,
note that each induced subgraph $G_K[V(G_{1, i})\cup \{b_i\}]$ (see
Fig.7(b)) contains at least two vertices belonging to $S$, it means
that $\gamma^{d}_2(G_K)\geq \frac{|G_K|}{3}$. Combining this
conclusion with Theorem~2.6, we have that
$\gamma^{d}_2(G_K)=\frac{|G_K|}{3}$. Since each graph of
$\mathscr{U}$ is a claw-free graph with minimum degree $2$, we have
that the upper bound of Theorem~2.6 can not be improved, even when
restricted to the claw-free graphs with minimum degree $2$.

%%%%%%%%%%%%%%%%%%%%%%%%%%%%%%%%%%%%%%%%%%%%%%%%%%%%%%%%%%%%%%%%%%%%%%%

\end{document}